\documentclass[12pt]{amsart}

\usepackage{amsfonts, amsmath, amssymb, verbatim}
\usepackage[mathcal]{euscript}

\newcommand{\Q}{\ensuremath{{\mathbb Q}}}
\newcommand{\Z}{\ensuremath{{\mathbb Z}}}
\newcommand{\R}{\ensuremath{{\mathbb R}}}
\newcommand{\N}{\ensuremath{{\mathbb N}}}
\newcommand{\C}{\ensuremath{{\mathbb C}}}

\newcommand{\XX}{\ensuremath{{\bold X}}}
\newcommand{\YY}{\ensuremath{{\bold Y}}}
\newcommand{\CC}{\ensuremath{{\bold C}}}

\newcommand{\Ass}{\ensuremath{{\bold {Ass}}}}

\newcommand{\AAA}{\ensuremath{{\bold A}}}

\newcommand{\WK}{\ensuremath{\XX^K}}

\newcommand{\Omm}{\ensuremath{\Omega(\WK)}}
\newcommand{\OmY}{\ensuremath{\Omega (\Sigma\YY)^K}}

\newcommand{\Hmm}{\ensuremath{H_*(\Omm)}}
\newcommand{\CCK}{\ensuremath{{\bold C}_K}}
\newcommand{\TCK}{\ensuremath{\widetilde{\bold C}_K}}

\newcommand{\PKN}{\ensuremath{{P_K^{\,0}}}}

\newcommand{\cat}[1]{\mbox{\sc #1}}
\newcommand{\Top}{{\cat{top}_*}}
\newcommand{\dga}{\cat{dga}}
\newcommand{\ga}{\cat{ga}}

\newcommand{\dg}{\dga}
\newcommand{\dgc}{{\cat{dgc}}}
\newcommand{\gc}{\cat{{gc}}}

\newcommand{\hopf}{{\cat{hopf}}}

\newcommand{\catk}{\cat {k}}

\newcommand{\colim}{\operatorname{colim}}

\newcommand{\supp}{\operatorname{supp}}
\newcommand{\skel}{\operatorname{skel}}
\newcommand{\pt}{\operatorname{{pt}}}
\newcommand{\Tor}{\operatorname{Tor}}
\newcommand{\Ext}{\operatorname{Ext}}
\newcommand{\Cotor}{\operatorname{Cotor}}
\newcommand{\Ker}{\operatorname{Ker}}
\newcommand{\id}{\operatorname{Id}}

\newtheorem{theorem}{Theorem}[section]
\newtheorem{proposition}[theorem]{Proposition}
\newtheorem{lemma}[theorem]{Lemma}
\newtheorem{corollary}[theorem]{Corollary}

\theoremstyle{definition}
\newtheorem{definition}[theorem]{Definition}
\newtheorem{construction}[theorem]{Construction}
\newtheorem{example}[theorem]{Example}
\newtheorem{examples}[theorem]{Examples}
\newtheorem{notation}[theorem]{Notation}

\theoremstyle{remark}
\newtheorem{remark}[theorem]{Remark}

\numberwithin{equation}{section}

\title{Loops on polyhedral products and diagonal arrangements}
\author{Natalia Dobrinskaya}
\address{VU University Amsterdam, De Boelelaan 1081, Amsterdam 1081HV, the Netherlands}
\email{NE.Dobrinskaya@few.vu.nl}

\begin{document}


\begin{abstract}
This is the first of a series of papers that investigates
the loop space homology of polyhedral products.

To any simplicial complex $K$ on $m$ vertices there corresponds a
polyhedral product functor, which associates to $m$ based
topological spaces ${\bf X}=(X_1,...,X_m)$ a certain subspace
$\WK$ in the cartesian product $\prod_i X_i$. In this paper we
establish a connection between the loop space homology of polyhedral
products of any 1-connected spaces and the homology of certain
diagonal arrangements associated with $K$. This reduces the
problem to the calculation of the $\Ext$-algebra of the exterior
Stanley-Reisner algebra of $K$. We illustrate these results by
finding the presentation of such loop homology algebras for flag
complexes and skeletons of simplices, generalizing results of
Panov-Ray, Papadima-Suciu, Lemaire.

Finally, we show that in the case when all the $X_i$'s are suspensions,
the homology splitting comes from the stable homotopy splitting of $\Omm$.
\end{abstract}

\maketitle

\renewcommand{\baselinestretch}{1.0}

\renewcommand{\baselinestretch}{1.3}

\section{Introduction}
Let K be a simplicial complex on the vertex set
$[m]=\{1,\dots,m\}$ and let $\XX=(X_1,\dots,X_m)$ be a sequence of
based topological spaces. The polyhedral product, or {\it
$K$-product}, is a natural subspace $$\WK=(X_1,\dots,X_m)^K\subset
X_1\times\dots \times X_m$$ defined by the following condition:
\begin{multline*}
(x_1,\dots,x_m)\in \WK \Leftrightarrow\\
\text{ for any } \tau\notin K \text{ there exists } i\in\tau
\text{ such that } x_i \text{ is the base point of } X_i.
\end{multline*}

This construction was introduced in \cite{anick} and
generalizes usual wedges (when $K$ is the set of vertices), cartesian
products ($K$ is $(m-1)$-dimensional simplex $\Delta^{m-1}$), fat
wedges ($K=\partial \Delta^{m-1}$). Our interest in it arises from
certain problems in toric topology where $K$-products play an
important role: for example, so called Davis-Januszkiewicz spaces $DJ_K$
are constructed in this way as $(\C P^\infty,\dots,\C
P^\infty)^K$.  In this paper we consider the problem of
calculating the loop homology algebra of this construction.

It is well known that taking homology of the loop space of a given
based topological space with field coefficients preserves products
and coproducts, which for topological spaces are wedges and
cartesian products; here homology is regarded as a graded Hopf
algebra with Pontryagin multiplication. It turns out that for
some strict class of simplicial complexes $K$ (so called flag
complexes) the answer is as simple as for wedges and products: the
required algebra is the colimit of the corresponding diagram in
the category of graded algebras. In general this answer fails due
to non-triviality of higher order Whitehead products in $\WK$.
In the present paper we investigate the functor for the loop space
homology of $K$-products.

Another interesting question raised in \cite{BP:book} concerns
possible relations of homology calculations for $\Omega DJ_K$ and
for the diagonal arrangement associated with $K$. In this paper we
explain this connection. We show that the diagonal arrangement
plays an important role in constructing geometric models for the
loop spaces for all $K$-products and, moreover, this leads to the
corresponding decomposition results in homology.

Our approach is inspired by the theory of labelled configuration
spaces. The results obtained are similar to the connection between
iterated loop spaces and classical configuration spaces such as
Milgram and May results, Snaith stable splitting.

We define the configuration space $\CCK=\sqcup_{I\in \N^m}\CCK(I)$
of particles in $\R^1$ with labels and partial collisions (see the
precise definition in Section~\ref{Sec:cktk}). This space appears
in the configurations space models for our loop spaces:
$$\Omm\simeq \sqcup \CCK(I)\times { \Omega {\bf X}}^I/\sim,$$
which we discuss in Section~\ref{Sec:split}. This model is
obtained using methods developed in \cite{dobrf} and originating
with Segal.

In the particular case when each space of $\XX$ is a suspension:
$X_i=\Sigma Y_i$, $i\in [m]$, we get a stable splitting (see
Theorem~\ref{T:split}):
$$\Omega (\Sigma Y_1,\dots,\Sigma Y_m)^K\simeq_{s}\bigvee _{I=(i_1,\dots,i_m)\in\N^m}
\CCK(I)_+\wedge Y_1^{\wedge i_1}\wedge\dots\wedge Y^{\wedge
i_m}.$$ which obviously leads to the homology splitting in this
case.

For general ${\bf X}$, we use the certain algebraic structures on
$H_*(\CCK)$ and our main homology result states

\begin{theorem}\label{T:main}
Let $X_1$,\dots,$X_m$ be based 1-connected topological spaces.
Then for homology with field coefficients the following graded
algebra isomorphism holds
\begin{equation}\label{Eq:thm}
H_*(\Omega \WK)\cong \bigoplus _{I\in\N^m}H_*(\CCK(I))\otimes
\tilde{H}_*(\Omega\XX)^{\otimes I}/\sim
\end{equation}
The equivalence relation is the closure of the relations
of the form
$$(\mu_{j,k}y)(\dots \otimes x^k_j\otimes x^{k+1}_j\otimes\dots)\sim y(\dots\otimes x^k_j\cdot x^{k+1}_j\otimes\dots),$$
where $y\in H_*(\CCK(I))$, $x_j^k,x_j^{k+1}\in H_*(\Omega X_j)$
and $\mu_{j,k}:\CCK(I)\to \CCK(I+e_j)$ are certain doubling
operations defined in Section~\ref{Subsec:doubling}.

If none of the algebras $H_*(\Omega\XX)$ has torsion, then the
isomorphism holds also over the integers.
\end{theorem}
The algebra structure on the right side is defined by the algebra structure on $H_*(\CCK)$ constructed in 
Section~\ref{Sec:tk} and the usual rule for tensor products.

This theorem looks more elegant when formulated in operadic
language:

 \noindent{\bf Theorem~1.1'} \ {\it If $\XX$ is as in
Theorem 1.1, then
$$H_*(\Omega\, \WK)\cong H_*(\CCK) \otimes _{\Ass} \tilde{H}_*(\Omega \XX),$$
where $\Ass$ is the non-$\Sigma$ associative operad, and
$H_*(\CCK)$ has the natural structure of a left module and a right
multi-module over $\Ass$.}

This operad description of the theorem will be discussed somewhere
else, and trying to make this paper more accessible, we avoid here
the operad language.

Due to Theorem~\ref{T:main}  the original problem is now reduced
to the calculation of the homology algebra of $\CCK$ and the
action of doubling operations $\mu$ on it.

The simple answer for flag complexes has the following explanation:
{\it for any flag complex $K$, all connected components of $\CCK$
are contractible}, and the presentation is this case is as follows
(here $\sqcup$ denotes the coproduct of connected graded algebras,
see Section~\ref{Sec:hocolim}):

\begin{corollary}\label{Cor:flag}
For any flag $K$ the following isomorphism holds for homology with
field coefficients
$$H_*(\Omega\, \WK)\cong
\sqcup_{i=1}^m H_*(\Omega X_i)/\sim,$$ with $x\cdot y\sim
(-1)^{\deg x\deg y}y\cdot x$ for $x\in H_*(\Omega X_i)$, $y\in
H_*(\Omega X_j)$ when $\{i,j\}\in K$.
\end{corollary}

This corollary was obtained rationally for $X=\CC P^\infty$ in
\cite{panray}, and for $X=S^{2n+1}$ in \cite{papasuc}.

For all non-flag $K$ this presentation fails to be true due to
higher homology of $\CCK$. For example, for each minimal
non-simplex $\tau$ of $K$ with no less than 3 vertices we have a
non-trivial $(|\tau|-2)$-dimensional class in $H_*(\CCK)$. This
class gives rise to the non-trivial $s$-product of elements from
$H_*(\Omega X_i)$, $i\in \tau$ and provides new generators in
$\Hmm$. This higher $s$-product in loop homology corresponds to
the higher order Samelson product in $\pi_*(\Omm)$ or, equivalently, to
the higher order Whitehead product in $\pi_*(\WK)$. We investigate the
homology of $\CCK$ in terms of higher operations in
\cite{Dobr:Poisson}, and apply
Theorem~\ref{T:main} to the discovered structures.

The connection to the theory of subspace arrangements is provided
by the fact that the components $\CCK(I)$ can also be expressed as
{\it complements of diagonal subspace arrangements} associated
with $K$. Thus, the problem of calculating the homology algebra
$\Hmm$ is now related to the problem of calculating the homology
of the complements of certain diagonal arrangements.

In Section~\ref{Sec:tck} we introduce a certain homotopy
modification $\TCK$ of the space $\CCK$ more suitable for homology
calculations. We show in Section~\ref{Sec:tk} the isomorphism of
its cellular chain complex with the Adams-Hilton model for the standard cell decomposition of
$(S^1,\dots,S^1)^K$ and with the cobar construction on
$H_*(S^1,\dots,S^1)^K$. This implies
\begin{theorem}\label{T:cotor}
For the configuration space $\CCK$ the following multi-graded
algebra isomorphism holds
$$\oplus_{n,I}H_{n,I}(\CCK; \Z)\cong \oplus_{n,I} \Ext_{|I|-n,I}^{\wedge(K)}(\Z;\Z).$$
\end{theorem}
Here $\wedge(K)$ denotes the exterior Stanley-Reisner algebra,
which is the factor of the exterior algebra
$\wedge[v_1,\dots,v_m]$ by the Stanley-Reisner ideal (see
Definition~\ref{Def:SR}).

In particular, the component $\CCK(1,\dots,1)$ corresponds to the
complement $D_K$ of the standard diagonal arrangement associated
with $K$, and this fact leads to the expected relation between
$H_*(D_K)$ and $H_*(\Omega DJ_K)$ (see
Section~\ref{Subsec:DJdiag}).

In this paper we illustrate an application of
Theorem~\ref{T:main} calculating the presentation of the algebra
$H_*(\Omega DJ_K)$ and the Lie algebra $\pi_*(DJ_K)$ for
$K=\skel_i\Delta[m]$ which is the $i$-skeleton of the full simplex
 (Proposition~\ref{Prop:skel}). We discuss the connection of
these calculations with known results about $(i+2)$-equal
arrangements (Corollary~\ref{Cor:equal}).

Many of the results of this paper were obtained during the
academic year 2006/2007, and the author is grateful to the
Universit\'e Catholique de Louvain for the hospitality, and
personally to Yves Felix and Pascal Lambrechts
for helpful discussions, references and
Pascal's corrections of the early draft.
I would also like to thank
Tilman Bauer, Dietrich Notbohm, Alvise Trevisan and Sadok Kallel
for discussions and corrections,
and especially Victor Turchin
for careful reading and useful suggestions.

\section{Preliminaries}

\subsection{Polyhedral products}\label{Sec:hocolim}
Here we give the definition of polyhedral product of based
topological spaces and outline some known results.

Let $K$ be an abstract simplicial complex, that is the collection
of subsets of $[m]==\{1,\dots,m\}$ which includes the empty set
and is closed under taking subsets. Let $\XX=(X_1,\dots,X_m)$ be a
sequence of based topological spaces.

Define $\XX^\sigma$ as $\prod_{i\in \sigma} X_i$ for each
$\sigma\in K$. Then the inclusion $\sigma\subset [m]$ let us
regard $\XX^\sigma$ as the subspace in the full product:
$$\XX^\sigma\subset  \prod_{i\in [m]}X_i.$$
Now the definition of $K$-product in the introduction is
equivalent to the following one.
\begin{definition}
The polyhedral product, or $K$-product, of $\XX=(X_1,\dots,X_m)$
is
$$\WK=\cup_{\sigma\in K}\XX^\sigma.$$
In the particular case when $X_1\cong\dots\cong X_m\cong X$, we call
$\WK$ the $K$-power of $X$ and denote it by $X^K$.
\end{definition}

\begin{examples}
\begin{enumerate}
\item {\bf Product.} If $K=\Delta[m]$ is the full $(m-1)$-dimensional simplex on $[m]$, then $$\WK\cong X_1\times\dots \times
X_m;$$
\item {\bf Wedge.} If $K$ is the set of disjoint vertices: $K=\{\emptyset,
\{1\},\dots,\{m\}\}$ then $$\WK\cong X_1\vee\dots
\vee X_m;$$
\item {\bf Fat wedge.} If $K=\partial \Delta[m]$ is the boundary of a $m-1$-dimensional simplex then $\WK$ is the fat wedge of $\XX$:
$$\WK\cong T_{m-1}(X_1,\dots,X_m);$$
\item {\bf Generalized fat wedge.} If $K=\skel_i \Delta[m]$ is the $i$-th skeleton of
the full simplex, then $\WK$ is a certain
generalization of the fat wedge (see e.g. \cite{porter}):
$$\WK\cong T_{i+1}(X_1,\dots,X_m);$$
\item {\bf Davis-Janiszkiewicz spaces.} If $K$ is arbitrary, and all $X_i\cong \CC P^\infty$, then $\WK$ is known as
Davis-Januskiewicz space $DJ_K$ and it plays an important role in
toric topology, being homotopy equivalent to the Borel construction
of a (quasi)-toric manifold $M^{2n}$:
$$(\C P^\infty)^K\simeq DJ_K\simeq  M^{2n}\times_{T^n} ET^n.$$
\end{enumerate}
\end{examples}

The definition can be reformulated in categorical language passing to the
following diagram. Let $\catk$ be the small category whose objects
are all the simplices of $K$ (including the empty simplex), and the
morphisms are the inclusions of simplices.

\begin{definition}
The exponential diagram $S_K(\XX)$ in $\Top$ is the functor
$\catk\to \Top$ which assigns $\XX^\sigma$ to each simplex
$\sigma$, and the natural inclusion
$\XX^\sigma\hookrightarrow\XX^\tau$ to each morphism
$\sigma\subset \tau$.
\end{definition}

Now it is easy to see that
$$\WK=\colim^\Top S_K(\XX).$$

 The analogous diagrams can be considered in different
categories which are augmented and have product and coproduct, for
example $\dgc$, $\hopf$. Throughout the paper we use the following
notation for categories.

\begin{notation}

\begin{itemize}
\item
$\Top$: based topological spaces;\\ coproduct is the wedge,
product is the cartesian product;
\item
$\cat{dgc}$: coaugmented differential graded coalgebras;\\
$\cat{gc}$: connected graded coalgebras,\\ coproduct is the
connected direct sum, product is the tensor product;
\item
$\hopf$: graded connected cocommutative Hopf algebras;\\
coproduct is the free product, product is the tensor product;
\item
$\ga$: graded connected algebras;\\  $\dga$: augmented
differential graded algebras;\\ coproduct is the free product, and we will not need the product as it is explained below.
\end{itemize}
\end{notation}

We will use the same definition of the diagram $S_{\bf K}$ in
these categories with one exception. Writing $S_K(H_*(\Omega
\XX))$ we use the tensor product which is the product in the
category of cocommutative Hopf alebras, but when we are interested
only in algebra structure on the colimit of this diagram, we write
$\colim^{\dga} S_K(H_*(\Omega \XX))$. This means that we just
forget the coalgebra structure, as coproduct in the category
$\dga$ is also the free product.

It is easy to check that
$$C_* (\WK)\simeq \colim^\dgc S_K(C_*(\XX)),$$ where $C_*:\Top\to \dgc$ is
the chain coalgebra functor. The examining differentials in
the corresponding Mayer-Vietoris sequence implies that over a
field the following isomorphism holds
$$H_* (\WK) \simeq \colim^{\gc}S_K(H_*(\XX))),$$
where the colimit is taken in the category $\gc$.

The last isomorphism can be also derrived from the stable homotopy splitting proved in \cite{cohenetc}:
$$\Sigma \WK\simeq \Sigma\bigvee_{\sigma\in K}\XX^{\wedge \sigma}.$$

It turns out that colimit of the diagrams $S_K$ does not commute
with the loop functor. We cannot see that in the simplest cases as
the following examples show. Here $R$ denotes a field.

\begin{examples}
\begin{enumerate}
\item For the cartesian
product ($K=\Delta[m]$) we have the following isomorphism
$$H_*(\Omega(X_1\times\dots \times X_m);R)\cong
H_*(\Omega X_1;R)\otimes\dots\otimes H_*(\Omega X_m;R)$$
\item \cite{CS} For the usual wedge ($K=\{\emptyset,
\{1\},\dots,\{m\}\}$) we have
$$H_*(\Omega(X_1\vee\dots \vee X_m);R)\cong H_*(\Omega
X_1;R)\sqcup\dots\sqcup H_*(\Omega X_m;R),$$ where $\sqcup$
denotes the coproduct for connected graded algebras (free
product).
\end{enumerate}
\end{examples}

The smallest complex for which the two functors do not commute is the
boundary of the 2-simplex.

\begin{example}
Let $K=\partial \Delta[3]$ be the boundary of the 2-dimensional
simplex, and $X_1=X_2=X_3=S^n$, $n\ge 2$. For any commutative ring
$R$ we have
$$
H_*(\Omega(\WK);R)\cong R[u_1,u_2,u_3]\sqcup R[w] \ncong
R[u_1,u_2,u_3]=\colim^{\dga} S_K(H_*(\Omega\XX);R),$$ where $\deg
u_i=n-1$, $\deg w=3n-2$. Here $R[v_1,\dots,v_k]$ denotes the
polynomial algebra generated by $v_1,\dots,v_k$.
\end{example}

This paper is therefore devoted to finding the correct functor computing  $H_*(\Omm)$ via
$H_*(\Omega\XX)$ in the category of graded algebras.

\subsection{}
Throughout the paper we will use the following notation.

\begin{notation}
\begin{itemize}
\item
$[n]$ is the set $\{1,2,\dots,n\}$;
\item
$\N$ is the semigroup of non-negative integers;
\item
$e_j=\{0\}^{j-1}\times \{1\}\times \{0\}^{m-j}\in \N^m$ for $j\in
[m]$;
\item
for $I=\{i_1,\dots,i_m\}\in \N^m$, $|I|=i_1+\dots+i_m$;
\item
for $I=\{i_1,\dots,i_m\}\in \N^m$, $\supp I=\{j\in[m]\,|\,i_j\ne
0\}$;
\item
for topological spaces $\XX=(X_1\dots,X_m)$ and
$I=(i_1,\dots,i_m)\in \N^m$
$$\XX^I=X_1^{\times i_1}\times\dots\times X_m^{\times i_m};$$
\item
for (graded) vector spaces ${\bf V}=(V_1,\dots,V_m)$ and
$I=(i_1,\dots,i_m)\in\N^m$
$${\bf V}^I=V_1^{\otimes i_1}\otimes\dots\otimes V_m^{\otimes i_m};$$
\item
for two tuples $I,J\in \N^m$ we write $I\prec J$ if $i_k\le j_k$
for every $k\in [m]$ and $I\ne J$;
\item
for a commutative ring $R$ the polynomial algebra over $R$
generated by $v_1,\dots,v_m$ is denoted by $R[v_1,\dots,v_m]$;
\item
for a commutative ring $R$ the exterior algebra over $R$ generated
by $v_1,\dots,v_m$ is denoted by $\wedge[v_1,\dots,v_m]$.
\end{itemize}
\end{notation}

\subsection{Combinatorial definitions and notation}
First define certain special simplicial complexes.
\begin{notation}
\begin{itemize}
\item
$\Delta[S]$ is the full simplex on the vertex set $S$;\\
$\Delta[m]$ denotes $\Delta[S]$ for $S=[m]$;
\item
$V[S]$ is the simplicial complex on $S$ consisting of a disjoint union of the vertices ;\\
$V[m]$ denotes $V[S]$ for $S=[m]$;
\item
$\skel_i\Delta[m]$ denotes the simplicial complex which consists
of all subsets $\sigma\subset[m]$ of cardinality no more than
$i+1$; geometrically it means that it is the collection of all
simplices from $\Delta[m]$ with the dimension $\le i$.
\end{itemize}
\end{notation}

Let further $K$ be an abstract simplicial complex on $[m]$.

\begin{notation}
\begin{itemize}
\item
$\bar{K}$ denotes the set of all non-empty simplices of a a
simplicial complex $K$.
\item
for a simplex $\sigma$ of a simplicial complex $K$, $|\sigma|$
denotes the number of elements in $\sigma$, and $\dim\sigma =
|\sigma|-1$ its dimension.
\item
for any subset $V\subset [m]$, $K|_V$ denotes the full subcomplex
$\{\sigma\in K|\sigma\subset V$\};

\end{itemize}
\end{notation}

We associate with each simplicial complex $K$ two algebras in teh
following way.

\begin{definition}\label{Def:SR}
The Stanley-Reisner algebra of a simplicial complex $K$ over a
commutative ring $R$ is defined as
$$R(K)=R[v_1,\dots,v_m]/I_{SR}(K),$$
where $I_{SR}(K)$ is the Stanley-Reisner ideal defined by
$$v_I:=v_{i_1}\dots v_{i_s}=0 \,\text{ if }\, I=\{i_1,\dots,i_s\}\notin K.$$

The exterior Stanley-Reisner algebra $\wedge(K)$ is defined
similarly by
$$\wedge(K)=\wedge[v_1,\dots,v_m]/I_{SR}(K).$$
\end{definition}

Note that in the above definition to construct theh ideal $I_{SR}$
it is enough to take only minimal non-simplices of $K$.

\begin{definition}
A subset $\tau\subset K$ is called a {\it missing face} if
$\tau\notin K$ but any proper subset $\sigma\varsubsetneq \tau$ is
in $K$. Missing faces are therefore the minimal non-simplices in
$K$. Geometrically it means that $\tau\notin K$ but $\partial
\tau\in K$.
\end{definition}

Now we are ready to define a certain subclass of simplicial
complexes, which turned out to be special for the problem
considered in the present article.

\begin{definition}\label{Def:flag}
$K$ is called {\it flag} if every its missing face has 2 elements or, equivalently, dimension 1.
\end{definition}
Obviously, a flag complex $K$ is completely determined by its
1-skeleton $\skel_1 K$.

\begin{examples}
\begin{enumerate}
\item The boundary of any polygon with no less than 4 vertices is flag;
\item the $n$-dimensional octahedron for any $n$ is flag;
\item the boundary of the $n$-simplex is non-flag for any $n$;
\item more generally, the $i$-dimensional skeleton of the $n$-simplex for $n>i\ge 1$ is non-flag.
\end{enumerate}
\end{examples}

In certain statements and proofs we use the following construction which assigns to
any simplicial complex with $s$ vertices the $s$-functor for
simplicial complexes.
\begin{construction}\label{Cons:functorK}
Suppose we have a simplicial complex $K$ on $[s]$, and $s$
simplicial complexes ${\bf L}=\{L_1,\dots,L_s\}$ on
$V_1$,\dots,$V_s$ respectively. This allows to define a functor
which we denote by $K(L_1,\dots,L_s)$, or $K({\bf L})$, which is a
simplicial complex on $V({\bf L})=V_1\sqcup\dots\sqcup V_s$
constructed by the following rule. Any $I\subset V({\bf L})$ is
naturally defined by the data: $\sigma(I)\subset [s]$ and the
non-empty sets $S_i(I)\subset V_i$ for each $i\in \sigma(I)$, such
that $I=\sqcup_{i\in \sigma(I)} S_i(I)$. Then the defining
condition is as follows: $I\in K(L_1,\dots,L_s)$ {\it if and only
if} $\sigma(I)\in K$ and each $S_i(I)\in L_i$.
\end{construction}

\section{Configuration space $\CCK$ and diagonal arrangements}\label{Sec:cktk}

\subsection {Construction of $\CCK$ as a configuration space}
In this section we construct the configuration space $\CCK$.
Informally speaking, this is the unordered configuration space of
particles on the real line, each labelled with a nonempty simplex from
$K$, and the topology is defined by the following rule: if two
particles with labels $\sigma$ and $\tau$ are getting close, they
can collide only in case when their labels $\sigma$ and $\tau$
are disjoint and $\sigma\sqcup\tau\in K$, and so the resulting
particle gets the label $\sigma\sqcup\tau$.

Also we can regard this space as the configuration space of
particles in $\R^1$ with labels from $[m]$, and the condition specifying which
subsets can have the same coordinate is determined
by the simplicial complex $K$. We give here the strict definition
using this second approach.

Let $B(k)$ denotes the classical unordered configuration space of
$k$ particles on the real line:
$$B(k) = \{{\bf t}=\{t(1),\dots,t(k)\}\subset \R\,|\, t(i)< t(j)\text{ for  }i< j\}.$$
 Construct the space $$\CCK=\sqcup_{I\in \N^m} \CCK(I),$$
where for each $I=\{i_1,\dots,i_m\}$ the component $\CCK(I)\subset
B(i_1)\times\dots\times B(i_m)$ is defined by the following
condition:
$$({\bf t}_1,\dots,{\bf t}_m)\in \CCK(I)\Leftrightarrow
 \cap_{j\in \tau} {\bf
t}_j=\emptyset \,\,\,\,{\rm for\,\, each}\,\,\,\, \tau\notin K.$$

We call a tuple $I$ the {\it multi-degree} of a configuration if
it is from $\CCK(I)$.
Note that each $\CCK(I)$ is not necessarily connected,
$\pi_0(\CCK)$ will be calculated in \ref{Sec:tck}.

\subsection{$\CCK$ as complements of diagonal
arrangements}\label{Subsec:ckdiag}

Given a simplicial complex $K$ on $[m]$ we associate with it the
arrangement of diagonal subspaces  in $\R^m$ as follows. It
consists of the subspaces $$x_{j_1}=x_{j_2}=\dots=x_{j_s},$$ one for each
missing face $j=\{j_1,\dots,j_s\}$ of $K$. The complement
of this arrangement we denote by $D_K$. The direct comparison of
the definitions shows that
$$D_K\cong \CCK(1,\dots,1).$$
The components of other multi-degrees also admit similar
descriptions, we summarize this in the following lemma.

\begin{proposition}\label{Prop:diag}
The components $\CCK(I)$ have the following description in terms of
diagonal arrangements:
\begin{enumerate}
\item
The component $\CCK(1,\dots,1)$ is homeomorphic to the complement
of the diagonal arrangement $D_K$ associated with $K$.
\item
The component $\CCK(I)$ for $I\prec (1,\dots,1)$ is homeomorphic
to the complement of the diagonal arrangement $D_L$ associated
with the full subcomplex $L=K|_{\supp(I)}$, where $\supp(i_1,\dots,i_m)=\{j|i_j\ne 0\}$.
\item
The component $\CCK(I)$ for $I\npreceq (1,\dots,1)$ is
homeomorphic to a certain collection of connected components
of $D_{K'}$ associated with the complex
$K'=K(V[i_1],\dots,V[i_m])$. More precisely, this collection is
defined by the diagonal inequalities:
$$t_j(1)<t_j(2)<\dots<t_j(i_j),\, j\in[m]$$
\end{enumerate}
\end{proposition}

\begin{examples}\label{Ex:diagarr}
\begin{enumerate}
\item
The diagonal hyperplane arrangement which consists of all
hyperplanes $x_i=x_j$ corresponds to the simplicial complex
$K=V[m]$ --- the disjoint union of vertices from $[m]$.
\item
The so called ''$k$-equal arrangement'', consisting of all subspaces
of the form: $$x_{i_1}=\dots=x_{i_k},$$ corresponds to
$(k-2)$-skeleton of the simplex on $[m]$. Thus its complement,
called ''no $k$-equal manifold'' is homeomorphic to
$C_{\skel_{k-2}\Delta[m]}(1,\dots,1)$.
\item
The complement of any coordinate subspace arrangement in $\R^m$ is
homotopy equivalent to $D_K$ for some $K$: we add one additional
coordinate $x_{m+1}$ and replace the equation
$x_{i_1}=\dots=x_{i_k}=0$ by $x_{i_1}=\dots=x_{i_k}=x_{m+1}$.
Hence, it is homotopy equivalent to the $(1,\dots,1)$-component of
the corresponding configuration space $\CCK$.

More precisely, if we denote by $Z_K$ the complements of the
coordinate subspace arrangement with the subspaces of the form
$$x_{j_1}=\dots=x_{j_s}=0,$$
one for each missing face $\{j_1,\dots,j_s\}$, then
$$Z_K\simeq C_{\Sigma K}(1,\dots,1),$$
where $\Sigma K$ denotes the simplicial complex on $[m+1]$ determined by the
following condition: $\tau'$ is a missing face for $\Sigma K$ if
and only if $\tau'=\tau\sqcup\{m+1\}$ for some missing face $\tau$
of $K$.
\end{enumerate}
\end{examples}

\subsection{Monoid $\TCK$}\label{Sec:tck}
Here we construct a cell complex that is homotopy equivalent to
$\CCK$ but has the structure of a monoid. Moreover, we will see in
Section~\ref{Sec:tk} that it has a nice cellular
structure.

\begin{definition}
Let $\TCK$ be the space of tuples $(c_0,t_1,c_1,\dots,t_k,c_k)$
for all integer $k\ge -1$ (for $k=-1$ we take an empty tuple),
where $c_i\in K$, $t_i\in [0,1]$, which satisfy the following
property:\begin{itemize} \item[]
 if $t_j<1$,
$t_{j+1}<1$,\dots,$t_{j+s}<1$ for some $j$, $s$ then
$c_j,c_{j+1},\dots,c_{j+s}$ are pairwise disjoint and
$$c_j\sqcup c_{j+1}\sqcup\dots\sqcup c_{j+s}\in K,$$
\end{itemize}
 and with
the following identification:\item
$$(c_0,\dots,c_j,0,c_{j+1},\dots,c_k)\sim (c_0,\dots,c_j\sqcup
c_{j+1},\dots,c_k).$$
\end{definition}

From a more general theory developed by the author in \cite{dobrf}
(see Corollary~3.7) the following proposition follows.
\begin{proposition}
The space $\TCK$ is homotopy equivalent to $\CCK$. It has a
natural monoid structure:
$$(c_0,t_1,\dots,t_k,c_k)\cdot (c_0',t_1'\dots,t_s',c_s')=(c_0,t_0,\dots,c_k,1,c_0',t_0',\dots,c_s') $$
and its classifying space is
$$B\TCK\simeq (S^1)^K.$$
\end{proposition}

It is well-known that monoids $M$ such that $\pi_0(M)$ is a group,
satisfy the property $M\simeq \Omega BM$. It obviously fails for
$\TCK$, and it is not homotopy equivalent to a loop space.
However, we will see in Section~\ref{Sec:tk} that its cellular
chain complex is isomorphic to the cobar construction on
$H_*W_K(S^1)^K$.

In the rest of the section we find the presentation of the monoid
$\pi_0(\TCK)$. For any $K$ consider the {\it right-angled Artin
monoid} defined by the following presentation:
$$A^+_K=<y_1,\dots,y_m\,|\, y_iy_j=y_jy_i \text{ for } \{i,j\}\in K>.$$
Obviously this monoid depends only on the 1-skeleton of $K$.

\begin{lemma}
The following isomorphism of monoids holds:
$$\pi_0(\TCK)\cong \pi_0(\CCK)\cong A_K^+.$$
\end{lemma}

\begin{proof}
The required homomorphism $\pi_0(\TCK)\to A_K^+$ is given by taking product of all labels from the left to the right, and easily can be checked to be an isomorphism.
\end{proof}

Thus, if there are any missing edges in $K$, the spaces $\CCK(I)$
can be non-connected. Otherwise, when $\skel_1(K)$ is full, we
have that the isomorphism $\pi_0(\CCK)\cong\Z_+^m$ is given by
taking the multi-degree of a configuration.

\section{Cellular chain algebra of $\CCK$ and doubling operations}\label{Sec:tk}
\subsection{Cellular chain algebra}
In this paragraph we show that the Adams-Hilton model for $(S^1)^K$,
which we denote by $T_K$, is isomorphic to the cellular chain
algebra $C_*(\TCK)$ for a certain cell decomposition of $\TCK$,
and thus, $H_*(T_K)$ and $H_*(\CCK)$ are isomorphic as algebras.
As $T_K$ coincides with the cobar construction for the coalgebra
$H_*((S^1)^K)$, we get the description of $H_*(\CCK)$ as Ext of
the exterior Stanley-Reisner algebra (Theorem~\ref{T:cotor}).

\begin{definition}\label{Def:epsilon}
Let $T_K$ be a free tensor algebra generated by all the simplices
$\sigma$ in $\bar{K}$ (the empty simplex corresponds to the unity
of the algebra), degree of a generator $\sigma$ is its dimension
as a simplex in $K$: $\deg \sigma=\dim \sigma =|\sigma|-1$. We denote the
operation of tensor product by $|$. The differential is
defined by the formula:
\begin{equation}\label{Eq:diff}
d\sigma=\sum_{\substack{(\sigma_1,\sigma_2):\\
\sigma=\sigma_1\sqcup\sigma_2,\sigma_i\ne\emptyset}}
(-1)^{\epsilon (\sigma_1,\sigma_2)+|\sigma_1|}(\sigma_1|\sigma_2),
\end{equation}
 where
$\epsilon(\sigma_1,\sigma_2)$ is the number of pairs $(i,j)$ with
$i\in\sigma_1$, $j\in\sigma_2$ such that $i>j$.
\end{definition}

\begin{remark} The $\dg$-algebra $T_K$ is isomorphic to the Adams-Hilton model
constructed from the standard cell decomposition of $(S^1)^K$.
\end{remark}

Now we show that this is the right model for $\CCK$.
\begin{proposition}\label{Prop:ccktk}
The algebra $C_*(\TCK)$ is isomorphic to $T_K$.
\end{proposition}
\begin{proof}
We will show that $\TCK$ has a cellular decomposition such that
the corresponding differential chain algebra is isomorphic to
$T_K$.

The definition of $\TCK$ gives a natural cubical decomposition
in the following way. Consider a sequence
$(\tilde{\sigma}_1,\dots,\tilde{\sigma_s})$ where each
$\tilde{\sigma}_i$ is some simplex $\sigma_i$ together with the
order of its vertices. We associate with it the cube in $\TCK$
which consists of all tuples $(c_0,t_1,c_1,\dots,t_k,c_k)$ defined
by the following condition: for certain $0\le j_1<\dots<j_{s-1}<
k$ we have
$$\begin{array}{ccccc}
c_0\sqcup\dots\sqcup c_{j_1}&=&\sigma_1& ; &t_{j_{1}+1}=1\\
\dots\\
c_{j_{i-1}+1}\sqcup\dots\sqcup c_{j_i}&=&\sigma_k& ; &t_{j_i+1}=1\\
\dots\\
c_{j_{s-1}+1}\sqcup\dots\sqcup c_{k}&=&\sigma_k& \\
\end{array}
$$
The union of the cubes which are associated with the same
underlying tuple of simplices $(\sigma_1,\dots,\sigma_s)$ for all
possible orders of their vertices forms a cell which corresponds
to the element $(\sigma_1|\dots|\sigma_s)$ in $T_K$.
Combinatorially it is the product of $s$ permutohedra of
dimensions $|\sigma_1|-1,\dots,|\sigma_s|-1$.

Now it is sufficient to check that the boundary operator on this
set of cells is given by the formula~\ref{Eq:diff} which was used
to define the differential in $T_K$. This can be shown by the
direct calculation the boundary operation on each cube, and then
taking their sum.
\end{proof}

Recall that The exterior Stanley-Reisner algebra $\wedge(K)$ is
defined as
$$\wedge(K)=\wedge[v_1,\dots,v_m]/I_{SR}(K),$$
where $I_{SR}(K)$ is the Stanley-Reisner ideal (see
Definition~\ref{Def:SR}).

The exterior Stanley-Reisner coalgebra $\wedge^\vee(K)$ is defined as the dual of $\wedge(K)$.
Its basis consists of the elements $v_\sigma$ for each $\sigma\in \bar{K}$, and
the co-multiplication is given by the formula
$$\Delta (v_\sigma)=\sum_{\substack{(\sigma_1,\sigma_2):\\ \sigma=\sigma_1\sqcup \sigma_2}} (-1)^{\epsilon(\sigma_1,\sigma_2)}
v_{\sigma_1}\otimes v_{\sigma_2},$$ where
$\epsilon(\sigma_1,\sigma_2)$ was defined in~Definition~\ref{Def:epsilon}.
This coalgebra has the multi-grading in $\N^m$.

\begin{remark}
Obviously, the coalgebra $H_*((S^1)^K;\Z)$ is isomorphic to the
exterior Stanley-Reisner colagebra $\wedge^\vee(K)$ over $\Z$.
\end{remark}

The following statement is the direct corollary from the
definitions of the algebra $T_K$ and the cobar construction.
\begin{lemma}
$\dg$-algebra $T_K$ is isomorphic to the cobar construction on the coalgebra $\wedge(K)$.
The isomorphism preserves multi-grading.
\end{lemma}

As the corollary of the results of this section we get the isomorphism
$$H_{k}(\CCK(I),\Z)\cong \Cotor_{|I|-k,I}^{\wedge^\vee(K)}(\Z,\Z).$$
As this isomorphism respects the algebra
structure, this gives the proof of Theorem~\ref{T:cotor}.
\subsection{Doubling operations $\mu$}\label{Subsec:doubling} In this
section we define doubling operations which are used in the
statement of Theorem~\ref{T:main}. We give geometric construction of them
for the configuration space $\CCK$ and the action on the chain
algebra $T_K$.

Roughly speaking these operations are replacing one of the
particles with two particles on a very small distance
$\varepsilon>0$ and with the same label (colour). It is more
convenient to give strict definitions replacing the space $\CCK$
by the homotopy equivalent version of labelled small intervals,
and to define the action of little interval operad, but here we
can proceed as follows. Let the space $\CCK^\varepsilon$ be the
subspace of $\CCK$ with the following condition: the distance of
any two particles with the same colour is not less than
$\varepsilon$. Fix one of the particles of a configuration of
multi-degree $I$. This is equivalent of choosing the pair of
integers: $(j, k)$ with $j\in [m]$, $k\in [i_j]$, which means that
we fix the $k$-th particle among all the particles with color $j$
taken from left to right.

Now in each configuration of multi-degree $I$ add one more
particle of the colour $j$ on the distance $\varepsilon/2$ to the
right from the fixed one. This defined the map:
$$\CCK^\varepsilon(I)\to \CCK^{\varepsilon/2}(I+e_j).$$
As for any $\varepsilon>0$ the space $\CCK^\varepsilon(I)$ is
homotopy equivalent to $\CCK(I)$, this defines the following map
up to homotopy:
$$\phi_{j,k}:\CCK(I)\to \CCK(I+e_j).$$

Let construct this operations on the level of chains in $T_K$. Fix
the multi-degree $I$ and the pair $(j,k)$ as above. Take a basis
element $(\sigma(1)| \dots|\sigma(s))\in T_K(I)$. Find the minimal
$n$ such that exactly $k$ simplices among $\sigma(1)$, \dots,
$\sigma(n)$ contains $\{j\}$.

Set
$$\mu^T_{j,k}(\sigma(1)\,|\, \dots\,|\,\sigma(s))=\sum_{\substack{(\tau_1,\tau_2):\\ \tau_1\sqcup\tau_2=\sigma(n)-\{j\}}}(-1)^{\epsilon(\tau_1,\tau_2)}(\sigma(1)\,|\, \dots\,|\,\{j\}\sqcup \tau_1\,|\,\{j\}\sqcup \tau_2\,|\,\dots\,|\,\sigma(s)).$$

By linearity this map defines the homomorphism of vector spaces
$$\mu_{j,k}^T: T_K(I)\to T_K(I+e_j),$$
and it's easy to see that it commutes with the differential.

The following statement will be not used in proofs of the results
in the present paper, but it gives some geometric feeling of
$\mu$-operations. We omit its proof here.
\begin{proposition}
The map $$\phi_{j,k}:\CCK(I)\to \CCK(I+e_j)$$ and the homomorphism
$$\mu_{j,k}^T: T_K(I)\to T_K(I+e_j)$$ constructed above induce up to sign the same homomorphism
in homology:
$$H_*(\CCK(I))\cong H_*(T_K(I))\to H_*(\CCK(I+e_j))\cong H_*(T_K(I+e_j))$$
\end{proposition}

We denote this homomorphism in homology by $\mu_{j,k}$.

\section{Applications}
\subsection{Flag complexes}
Recall that flag complexes were defined in
Definition~\ref{Def:flag}.

Nice properties of loops of polyhedral products for flag complexes
were noticed by several authors. We give below the examples of
that. Here $T(U)$ denotes the free tensor algebra generated by
elements of $U$.

\begin{examples}
If $K$ is flag, then
\begin{enumerate}
\item(\cite{panray})
$H_*(\Omega (\C P^\infty)^K, \Q)\cong T(U)/(u_i^2=0,
u_iu_j+u_ju_i=0$ for $\{i,j\}\in K$), where $U=\{u_1,\dots,u_m\}$
with $\deg u_i=1$.
\item(\cite{papasuc})
$H_*(\Omega (S^{2k+1})^K, \Q)\cong T(U)/(u_iu_j-u_ju_i=0$ for
$\{i,j\}\in K)$, where $U=\{u_1,\dots,u_m\}$ with $\deg u_i=2k$.
\end{enumerate}
\end{examples}

Here we show that Theorem~\ref{T:main} implies that these results hold in
more general situation.

\begin{corollary}
Over any field
$$H_*(\Omega \WK)\cong
\sqcup_{i=1}^m H_*(\Omega X_i)/\sim,$$ with $$x\cdot y\sim
(-1)^{\deg x\deg y}y\cdot x$$ for $x\in H_*(\Omega X_i)$, $y\in
H_*(\Omega X_j)$ when $\{i,j\}\in K$.
\end{corollary}

In other words, for flag $K$ $$H_*(\Omega \WK)\cong\colim^{\dga}
S_K(H_*(\Omega\XX)).$$

The proof is direct from the following lemma.

\begin{lemma}
For any flag simplicial complex $K$ each connected component of
$\CCK$ is contractible.
\end{lemma}
\begin{proof}

The statement can be obtained by induction on number of missing
edges. We use the description of $\CCK$ in terms of diagonal
arrangements stated in \ref{Prop:diag}.

In the beginning we have the contractible open cone in $\R^{|I|}$,
and on each step we will get the union of the interior of disjoint
cones which are bounded by hyperplanes passing through 0. Adding
one more missing edge, we sect the complement by the corresponding
hyperplane. Any cone in the complement either doesn't intersect
this hyperplane, or it is divided into two new cones, and both of
them are again contractible.
\end{proof}

\subsection{Davis-Januszkiewicz spaces and diagonal
arrangements}\label{Subsec:DJdiag}
 The important examples of
$K$-product are provided by so called Davis-Januszkiewicz spaces
which is $K$-power of $\C P^\infty$:
$$DJ_K:=(\C P^\infty)^K.$$
The spaces $DJ_K$ arises in toric topology as Borel construction
for (quasi)-toric manifolds.

In the book~\cite{BP:book} Buchstaber and Panov observed that
there is some coincidence in cohomology calculations for loops on
Davis-Januszkiewicz spaces and for complements of real diagonal
arrangements.

Recall that $R(K)$ denotes the Stanley-Reisner algebra of the
simplicial complex $K$ over a ring $R$:
$$R(K)=R[v_1,\dots,v_m]/I_{SR},$$
where $I_{SR}$ is the Stanley-Reisner ideal (see Definition~\ref
{Def:SR}).

If we set the degree of each generator $v_i$ equal to 2, then the
following isomorphism holds for the loops homology of
Davis-Januszkiewicz spaces with field coefficients $R$ (see
\cite{BP:book}):
$$H^* \Omega DJ_K\cong \Tor^*_{R(K)}(R,R).$$

From the other side, if the degree of each $v_i$ is set to 1, then
we have the following isomorphism for cohomology with field
coefficients of the complement of the diagonal arrangement
associated with $K$ (see \cite{PRW}):
\begin{equation}\label{Eq:peeva}
H^* (D_K) \cong \Tor_{R(K)}^{m-*}(R,R)_{(1,\dots,1)}.
\end{equation}

Now using the Theorem~\ref{T:main} we can explain this accident.

By Theorem~\ref{T:cotor} we have the multigraded algebra
isomorphism $H_{*,\bar{*}}(\CCK;\Z)\cong\Cotor_{|{\bar
*}|-*,\bar{*}}^{\wedge^\vee(K)}(\Z,\Z)$, where the degrees of the
generators in $\wedge(K)$ are equal to 1.

As $\Omega \C P^\infty\simeq S^1$ has torsion-free homology,
applying Theorem~\ref{T:main} to $DJ_K=(\C P^\infty)^K$ we get
$$H_*(\Omega DJ_K)\cong \Cotor_*^{\wedge^\vee(K)}(\Z,\Z)\otimes_{\{\mu_{j,k}\}} \wedge[u_1,\dots,u_m],$$
where $\deg u_i=1$, $\otimes$ is multidegree-wise tensor product
over the action of all the operations $\mu_{j,k}$, and $\Cotor$ is
regraded as described above.

Now if we use the integral formality of $DJ_K$ proved in
\cite{NR},
and change the multi-grading of $\Z(K)$
setting $\deg v_i=1$, we get the following statement.

\begin{corollary}
There exists the natural epimorphism of multi-graded vector spaces
$$f: A=\Ext_{\wedge(K)}^{*,\bar{*}}(\Z,\Z)\to B=\Ext_{\Z(K)}^{*,\bar{*}}(\Z,\Z)$$
such that
\begin{enumerate}
\item
$$f(y\cdot y')=(-1)^\epsilon f(y)\cdot f(y'),$$ where $y\in A_{n,
I}$, $y'\in A_{n', I'}$ and $\epsilon=
\sum_{j\le k}i_ji_k'-n'\sum_ji_j$.
\item
$$\Ker f|_{A_{n, I}}=\bigoplus_{j,k} \mu_{j,k} (A_{n, I-e_j}).$$
\end{enumerate}
\end{corollary}

\begin{remark}
Let $\Z_\tau(K)=\Z(K)/(v_1^2,\dots,v_m^2)$ be the truncated Stanley-Reisner algebra.
Then the above epimorhism $f$ is the composition of two homomorphisms
$$\Ext_{\wedge(K)}^{*,\bar{*}}(\Z,\Z)\to \Ext_{\Z_\tau(K)}^{*,\bar{*}}(\Z,\Z)\to \Ext_{\Z(K)}^{*,\bar{*}}(\Z,\Z),$$
where the first one is the additive isomorphism satisfying (1), and the second one is an algebra epimorphism with the kernel (2).
These maps are constructed by applying Theorem~\ref{T:main} to the space $(S^{2n})^K$ with $n\ge 1  $.
\end{remark}

\begin{examples}
The simplest examples of finding the kernel of the above
epimorphism are as follows.
\begin{enumerate}
\item
for any generator $x_j$ of multi-degree $e_j$ in
$\Ext_{\wedge(K)}^{*,\bar{*}}(\Z,\Z)$ we have
$\mu_{1,1}(x_j)=x_j^2$. So this element should be in the kernel
of $f$, and we get the additional relation in
$\Ext_{\Z(K)}^{*,\bar{*}}(\Z,\Z)$: $$x_j^2=0.$$
\item
for any missing face $\tau$ in $K$ we have a generator
$\omega_\tau$ of multi-degree $\sum_{j\in \tau} e_j$ in
$\Ext_{\wedge(K)}^{*,\bar{*}}(\Z,\Z)$. As for $j\in \tau$ we have
$\mu_{j,1}(\omega_\tau)=[x_j,\omega_\tau]$,
we get that these commutators are also in the kernel of $f$, so
the new relations are
$$[x_j,\omega_J]=0,$$
for $j\in J$.
\end{enumerate}
\end{examples}

\begin{remark}
One of the examples of explicit calculations of those two
$\Ext$-algebras and the epimorphism of the above Corollary  will
be given in proofs of Lemma~\ref{Lem:lemaire} and
Proposition~\ref{Prop:skel} (with topological grading of the
generators).
\end{remark}

Thus, we see that $\Cotor^{\Z^\vee(K)}(\Z,\Z)$ is obtained from
$\Cotor^{\wedge^\vee(K)}(\Z,\Z)$ by rescaling and then adding some
additional relations. We do not have any relations in multi-degree
$(1,\dots,1)$ as for the image of any $\mu_{j,k}$-operation we
have $i_j\ge 2$. Thus, we have only the rescaling of homology
groups $H_*(DJ_K;\Z)_{(2,\dots,2)}$ and
$H_*(\CCK;\Z)_{(1,\dots,1)}\cong H_*(D_K;\Z)$, and this explain
the connection we were interested to find.

\subsection{Torsion in $\Hmm$}
For general $K$ the homology of the space $\CCK$ is not
torsion-free. This means that even in case when all  $H_*(\Omega
X_i;\Z)$ are torsion free, so that Theorem 1.1 is valid for
integer coefficients, the resulting homology $H_*(\Omm;\Z)$ can
have arbitrary torsion.

To show that we will use Example~\ref{Ex:diagarr}(3). Consider the
compliment $Z_L$ of the coordinate subspace arrangement associated
with a simplicial complex $L$ on $[m]$.
It was shown in \cite{bjorner2}  that the integral cohomology of
its complement can have arbitrary torsion.

Due to Example~\ref{Ex:diagarr} the complement of a coordinate
subspace arrangement is homotopy equivalent to $D_K$ for the
certain $K$. Thus, for that $K$, $H_*(\CCK(1,\dots,1))$
also can have any
torsion.

Nevertheless, for certain classes of simplicial complexes the
space $\CCK$ will be torsion free. The simplicial complex is
called {\it shifted} if there is an ordering on its vertices, say
$1<\dots<m$, such that whenever $\sigma\in K$, $j\in \sigma$ and
$i<j$, we have $(\sigma-\{j\})\cup\{i\}\in K$. The skeletons
$K=\skel_i\Delta[m]$ are examples of shifted complexes.

\begin{proposition}\label{Prop:shifted}
For shifted simplicial complexes $K$ the homology $H_*(\CCK;\Z)$
is torsion-free. Thus, if all $H_*(\Omega X_i;\Z)$ are
torsion-free, then $\Hmm$ also has no torsion.
\end{proposition}

\begin{proof}
The proposition can be proved directly, but we use here another
way relying on the homotopy result from \cite{GT}.

Consider the space $Y^n_K=(S^{2n+1})^K$ for $n\ge 1$. Denote by
$F^n_K$ the homotopy fiber of the embedding $Y^n_K\hookrightarrow
(S^{2n+1})^m$.

Due to Theorem~9.4 from \cite{GT} in case when $K$ is shifted,
this space is homotopy equivalent to the following wedge of the
form $F^n_K\simeq \vee_{i=1}^s \Sigma^{k_i}(\wedge_{j\in J_i}
\Omega X_j)$ for certain subsets $J_i\subset [m]$. We have the
splitting of the loop spaces $\Omega Y^n_K\simeq (\Omega
S^{2n+1})^m\times \Omega F^n_K$, which implies the isomorphism of
the vector graded space
$$H_*(\Omega Y^n_K;\Z)\cong \Z[u_1,\dots,u_m]\otimes T(V),$$
where $T$ denotes the free tensor algebra over $\Z$, and $V$ is
some set of generators. Thus, $H_*(\Omega Y^n_K, \Z)$ has no
torsion. Due to Theorem~\ref{T:main} over integers,  all algebras
$H_*(\Omega Y^n_K, \Z)$ are certain rescalings of the algebra
$H_*(\CCK,\Z)$. So we get that $H_*(\CCK,\Z)$ has no torsion as
well.
\end{proof}

\subsection{Some calculations}
In this subsection we illustrate how to apply Theorem~\ref{T:main}
for finding the presentations of $\Hmm$.

Suppose we know the presentation of $H_*(\CCK)$. To derive the
presentation for $\Hmm$ we need to investigate the action of
$\mu$-operations on $H_*(\CCK)$. The easy but helpful observation
is that it's enough to find this action only on generators.

We give the example of such arguments getting results about
rational homotopy groups and the homology of $\Omega DJ_K$, when
$K$ is the skeleton of a simplex. These calculations are
generalized in \cite{Dobr:Poisson} for more general class of
simplicial complexes, but for this particular case we rely on the
presentation of $H_*(T_K)\cong H_*(\CCK)$ derived
from~\cite{lemaire}:

\begin{lemma}\label{Lem:lemaire}
Let $K=\skel_{s-2}\Delta[m]$, $s\ge 3$. The algebra $H_*(T_K;\Z)$
has the following presentation:
\begin{itemize}
\item $m$ generators $x_j$, $j\in [m]$ of degree $0$;
\item $C_m^s$ generators $\omega_J$, taken for all $J\subset [m]$ with $|J|=s$, of degree $s-2$.
\item relations $x_ix_j=x_jx_i$ for $i,j\in [m]$.
\item relations
$$\sum_{j\in S}(-1)^{\epsilon(j,S)}[x_j,\omega_{S-\{j\}}]=0$$
for $S\subset [m]$ with $|S|=s+1$.
\end{itemize}
\end{lemma}
\begin{proof}
Over a field the presentation was obtained in \cite{lemaire}. But
$K=\skel_{s-2}\Delta[m]$ is a shifted simplicial complex, so due
to Proposition~\ref{Prop:shifted} $H_*(\CCK,\Z)$ has no torsion.
Moreover, all the generators of the presentation are integral
classes: $\omega_J$ corresponds in $T_K$ to the differential $dJ$
of the non-existing element $J\notin T_K$.

\end{proof}

Now we need to find action of operations $\mu$ on the generators.
As each element $w_J$ is nothing else as $dJ$ in the model $T_K$ constructed in Section~\ref{Sec:tk}, the direct verification shows that
for $j\in J$ we have
$$\mu_{j,k}(\omega_J)=[x_j,\omega_J].$$
As the algebra $H_*(\Omega \C P^\infty)\cong \wedge[u]$ ($\deg u=1$) has no torsion, we obtain the following presentation for integral homology algebra $H_*(\Omega DJ_K)$.

\begin{proposition}\label{Prop:skel}
Let $K$ be $(s-2)$-dimensional skeleton of the full simplex
$\Delta[m]$, $s\ge 3$. Then
$$H_*(\Omega DJ_K; \Z)\cong T(u_j,w_J)/{\mathcal R},$$
where $w_J$ are taken for all $J\subset [m]$ with $|J|=s$, and
$u_j$ are taken for all $j\in [m]$. The degrees are  $\deg
w_J=2s-2$, $\deg u_i=1$, and the set of relations $\mathcal R$ is
as follows:
\begin{enumerate}
\item
$[u_i,u_j]=0$ for $i,j\in [m]$.
\item
$[u_j,w_J]=0$ for any $j\in J$;
\item
for any $S\subset [m]$ with $|S|=s+1$
$$\sum_{j\in S}[u_j,w_{S-\{j\}}]=0.$$
\end{enumerate}
\end{proposition}

This immediately implies the presentation of the rational homotopy groups for $DJ_K$:

\begin{corollary}
For $K=\skel _{s-2}\Delta[m]$ we have the following presentation
for rational homotopy algebra of $DJ_K$:
$$\pi_*(DJ_K)\otimes\Q\cong {\mathcal Lie}(\tilde{u}_j,\tilde{w}_J)/\tilde{\mathcal {R}}$$
with
\begin{itemize}
\item
each classes $\tilde{u}_i$ is the image of the generator $\pi_2(\C
P^\infty)$ under the inclusion of $i$-th copy $\C
P^\infty\hookrightarrow DJ_K$, and so $\deg \tilde{u}_i=2$;
\item
for $J=\{j_1,\dots,j_s\}$ the homotopy class
$\tilde{w}_J=\{u_{j_1},\dots,u_{j_s}\}$, which corresponds to the rational higher
order Whitehead product of the classes $u_{j_1},\dots,u_{j_s}$, and so
$\deg \tilde{w}_J=2s-1$;
\item
the relations $\tilde{\mathcal {R}}$ are the same as relations ${\mathcal {R}}$ from
Proposition~\ref{Prop:skel} where $[\cdot,\cdot]$ denotes now the
usual Whitehead product.
\end{itemize}
\end{corollary}

\begin{proof}
All generators obtained in Proposition~\ref{Prop:skel} are the
images of Hurewicz homomorphism. Moreover, for any space $Y$ the
commutator of Pontryagin product on $H_*(\Omega Y)$ corresponds to
usual Whitehead product on $\pi_*(Y)$. It is straightforward that the classes $\tilde{w}_J$ are the modified higher order Whitehead 
products, which are defined similarly to Massey products in cohomology.
It is known that they coincide with the classical ones up to rational coefficients and the indeterminacy.
\end{proof}

The interesting corollaries can be obtained by restricting the calculations to the
homology of the complements of the corresponding diagonal
arrangements. In this particular case the diagonal arrangements are called $s$-equal
arrangements, and their complements --- no-$s$-equal manifolds. We
recover the following result about the homology of these
manifolds proved in \cite{baryshnikov} (see also \cite{bjorner}).

\begin{corollary}\label{Cor:equal}
The integral homology of the no-$s$-equal manifold in $\R^m$,
which is the complement of $s$-equal arrangement, is torsion free
and has the following basis: each element is coded by the sequence
of pairwise disjoint subsets of $[m]$:
$$(I_1, J_1, I_2, J_2,\dots, J_k, I_{k+1})$$
with the following properties
\begin{itemize}
\item
$I_i, J_i\subset [m]$;
\item
$|J_i|=s$ for $i\in [k]$;
\item
$(\sqcup_i I_i)\sqcup (\sqcup_i J_i)=[m]$;
\item
$\max \{j\in J_i\}<\max \{j\in I_{i+1}\}$ for each $i\in[k]$.
\end{itemize}
This element has the degree $(s-2)k$, and is geometrically
represented by the product of spheres given by the following
system of equations in $\R^m=\{(x_1,\dots,x_m)\}$:
\[
\begin{cases}
x_j=3mn+j \text{ for each } n\in[k+1] \text{ and } j\in I_n;\\
\sum_{j\in J_i}|x_j-(3n+2)m|^2=1 \text{ and } \sum_{j\in
J_i}x_j=(3n+2)ms \text{ for }n\in [k].
\end{cases}
\]
\end{corollary}
\begin{proof}
First, by Proposition~\ref{Prop:shifted} $H_*(\CCK)$ has no
torsion. Then its presentation is obtained by combining
Lemma~\ref{Lem:lemaire}, Proposition~\ref{Prop:diag} and
Proposition~\ref{Prop:ccktk}. The geometrical description follows
from the representing the elements $w_J$ as the unit sphere in
$\R^J\backslash \{0\}$ intersected with the hyperplane $\sum
x_i=0$.
\end{proof}

\subsection{Poincare series of $\Hmm$}\label{Sec:poin}

For example, from Corollary~\ref{Cor:flag} it follows that for
flag complexes the formulas expressing $P_{\Omm}$ via $P_{\Omega
X_i}$ have the same form.

\begin{corollary}
For flag $K$
$$P^{-1}_{\Omm}(t) = \sum_{\sigma\in K}\prod_{j\in\sigma}(P^{-1}_{\Omega X_j}(t)-1)$$
\end{corollary}

It turns out that general formula is in some sense deformation of
these "flag formulas". The coefficients of those deformations
depends only on $K$, and could be calculated by means of
commutative algebra. We will discuss the resulting formulas in the
forthcoming paper, giving here just one example of such
calculations.

\begin{example}
If $K$ is the $(s-2)$-skeleton of a simplex:
$K=\skel_{s-2}\Delta[m]$ the space $\XX^K$ is a generalized fat
wedge $T_{s-1}(X_1,\dots ,X_m)$. For $m\ge 3$ and $s\ge 1$ $K$
fails to be flag, and the formula is
$$P^{-1}_{\Omega T_{s-1}(X_1,\dots ,
X_m)}(t)=P^{-1}_{\Omega X_1}(t)\cdot\dots\cdot P^{-1}_{\Omega
X_m}(t)+ (-t)^{s-1}\sum_{\sigma:|\sigma|\ge s-1}\prod_{i\in
\sigma} (P^{-1}_{\Omega X_i}-1).$$
\end{example}

\section{Labelled configuration spaces and stable splittings of the loop
spaces}\label{Sec:split}
The construction of labelled
configuration spaces with collisions considered in this section is
a particular case of the theory of configuration spaces with
labels in partial monoid developed, for example, in \cite{dobrf}.
We will use this construction only for partial monoids obtained as
$K$-products.

\begin{definition}\label{Def:cy}
Let $\YY=(Y_1,\dots,Y_m)$ be a sequence of well-pointed
topological spaces. Denote by $C(\R^1, ({\bf Y})^K)$ the space
$$C(\R^1, {\bf Y}^K):=\sqcup \CCK(I)\times {\bf Y}^I/\sim,$$
where the equivalence relations are defined by the following
condition: if some label $y\in Y_i$ is a basepoint $y=\star$, then
the point can be removed from the configuration.
\end{definition}

This configuration space $C(\R^1, {\bf Y}^K)$ can be viewed as the
space of unordered configurations in $\R^1$ with labels from ${\bf
Y}^K$ and collisions defined by the following rule: the points
with labels $(\sigma, {\bf y})$ and $(\tau, {\bf z})$ can collide
if and only if the simplices $\sigma$ and $\tau$ are disjoint in
$K$, and $\sigma\sqcup\tau\in K$.

\begin{remark}
Obviously we have that $\CCK\cong C(\R^1,(S^0)^K)$, where $S^0$ is
a 0-dimensional sphere.
\end{remark}

The similar definition can be given for the sequences of monoids.

\begin{definition}\label{Def:c1}
Let ${\bf A}=(A_1,\dots,A_m)$ be a sequence of topological
monoids. Denote by $\hat{C}(\R^1, {\bf A}^K)$ the space which as a
set of points coincides with the result of applying
Definition~\ref{Def:cy} to the sequence of topological spaces
${\bf A}$ with the base-points --- identities of monoids, but with
different topology: the points with labels $(\sigma, {\bf a})$ and
$(\tau, {\bf b})$ can collide if and only if the
$\sigma\cup\tau\in K$, and the result of their collision is a
point with the label $(\sigma\cup\tau, {\bf a}\cdot {\bf b})$,
where ${\bf a}\cdot {\bf b}$ is coordinate-wise product.
\end{definition}

The fundamental result in the theory of labelled configuration
spaces is Segal's theorem (\cite{segal}). It was generalized to
the case of labelled configuration spaces with collisions, see
e.g. Theorem 3.3 in \cite{dobrf}. Here will give the corollaries
of the last theorem for our cases.

\begin{corollary}\label{Cor:wktop}
For a sequence of connected topological spaces
$\YY=(Y_1,\dots,Y_m)$ the following homotopy equivalence holds
$$C(\R^1,\YY^K)\simeq \Omega (\Sigma\YY)^K,$$
where $\Sigma\YY$ denoted the sequence $(\Sigma Y_1,\dots,\Sigma
Y_m)$.
\end{corollary}

\begin{corollary}\label{Cor:wkmon}
For a sequence of connected topological monoids ${\bf
A}=(A_1,\dots,A_m)$ the following homotopy equivalence holds
$$\hat{C}(\R^1,{\bf A}^K)\simeq \Omega (B{\bf A})^K,$$
where $B{\bf A}$ denoted the sequence of classifying spaces $(B
A_1,\dots,B A_m)$.
\end{corollary}

The last corollary can be rewritten as the following homotopy
equivalence: for 1-connected $\XX=(X_1,\dots,X_m)$ we have
$$\Omm\simeq \hat{C}(\R^1,(\Omega \XX)^K).$$

This fact is related to the isomorphism in Theorem~1.1, as
$C(\R^1,(\Omega \XX)^K)$ is constructed from the space $\CCK$ and
the monoids $\Omega \XX$ (see Definition~\ref{Def:c1}).

Now let restrict ourself to the second case, when each $X_i$ from
$\XX$ is homotopy equivalent to the suspension $\Sigma Y_i$ for
some sequence of $\YY=(Y_1,\dots,Y_m)$ of connected spaces, and so
Corollary~\ref{Cor:wktop} is valid.

Now we can prove the result about stable splittings of such loop
spaces.

\begin{theorem}\label{T:split}
Let each $X_i\simeq\Sigma Y_i$ for some connected $Y_i$, $i\in
[m]$. Then there is a stable homotopy equivalence
$$\OmY\simeq_s \bigvee_{I\in \N^m} \CCK(I)_+\wedge \YY^{\wedge I}$$
\end{theorem}
\begin{proof}
We adapt arguments from the proof in \cite{bod} of the classical
Snaith splitting.

Consider the configuration space $F=C(\R^1, (\YY)^K)$ with
collisions. By Corollary~\ref{Cor:wktop} we have that $D\simeq
\Omega (\Sigma\YY)^K$.

The space $F$ is constructed as $F=\cup_{I\in \N^m} F_I$ for
$F_I=\CCK(I)\times \YY^I/\sim$, where equivalence relations are
defined in~\ref{Def:cy}.

Construct the new space $\hat{F}$ by imposing on $\sqcup
\CCK(I)\times {\bf Y}^I/\sim$ the stronger equivalence relations:
if some label $y\in Y_i$ is a basepoint $y=\star$, then the
configuration is equivalent to the empty one, or, in other words,
the {\it whole} configuration disappears.

It has natural multi-grading in $\N^m$ and it splits as a wedge
$\hat{F}=\vee_{I\in\N^m} \hat{F}_I.$

There are natural projections $\pi: F_I\to \hat{F}_I$. The space
$\hat{F}_I$ is in fact the filtration quotient for $F_I$:
$$\hat{F}_I\cong F_I/F_{<I},$$
where $F_{<I}=\cup_{J: J\prec I}F_J$.

Our goal is to prove the stable homotopy equivalence:
$\Sigma^\infty F\simeq \Sigma^\infty \hat{F}$. First, following
\cite{bod} construct the map $F\to
\Omega^\infty\Sigma^\infty\hat{F}$. As $\hat{F}$ is connected the
space $\Omega^\infty\Sigma^\infty\hat{F}$ is homotopy equivalent
to the labelled configuration space $C(\R^\infty, \hat{F})$, so in
fact we will define a map
$$P:F\to C(\R^\infty, \hat{F}).$$
We will write the points of $C(\R^\infty, \hat{F})$ in a form
$\sum (p, \hat{f})$ for $p\in \R^\infty$ and $\hat{f}\in \hat{F}$.

Let $f\in F_I$, so $f$ can be written as $(c,{\bf y})$, where
$c\in \CCK(I)$ and ${\bf y}\in \YY^I$. Take some subconfiguration
$\alpha \subset c$ of multi-degree $J$ (obviously $J\prec I$). As
each $\CCK(J)$ is a submanifold of $\R^{j_1+\dots+j_m}$, fix an
embedding of $\CCK=\sqcup_{J\in \N^m} \CCK(J)$ to $\R^\infty$, and
let $\bar{\alpha}$ be an image under this embedding. Then we
define
$$P(c,{\bf b})=\sum_{\alpha\in c}(\bar{\alpha}, \pi(c,{\bf b})).$$
It could be easily checked that this map agrees with the
equivalence relations in the construction of $F$.

The rest of the proof which shows that this map is homotopy
equivalence repeats the same part in the proof of \cite{bod} for
the following filtration of $F$ and $\hat{F}$:
$$
F_k=\bigcup_{|I|\le k} F_I\text{ and  } \hat{F}_k=\bigvee_{|I|\le
k}\hat{F}_I.
$$
\end{proof}

\section{Proof of Theorem~\ref{T:main}}
\subsection{Scheme of the proof}
Let ${\bf A}=(A_1,\dots,A_m)$ be a sequence of $\dg$-algebra's.
The object of out investigation now is the algebra
$$\widetilde{W}_K({\bf A}):=\Omega_*\colim^{\dgc}S_K(B_*\AAA),$$ where $B_*:\dga\to \dgc$
denotes the bar construction, and $\Omega_*: \dgc\to \dga$ --- the
cobar construction, and the diagram $S_K$ is defined in
Section~\ref{Sec:hocolim}.

For the explicit construction for it we need the following
construction of the $\dg$-algebra $P_K$. Take the $K$-product
$\N^K\subset \N^m$ (the role of the base point is played by zero
in $\N$), and consider a tensor algebra on $\N^K$ (again the
symbol for tensor product is $|$). Define the differential by its
value on a tuple $y\in \N^K$ using the formula:
$$d(y)=\sum_{y=y_1+y_2,y_i\ne 0}
(-1)^{\epsilon(y_1,y_2)}\,(y_1| y_2),$$ where the decomposition
into a sum is induced by the additive structure in $\N^m$. The
degree is defined by $\deg y=|y|-1$. The resulting $\dg$-algebra
we will denote by $P_K$.

$P_K$ can be endowed with multi-grading in the following way. One
grading is in $\N^m$, and it is calculated by taking the sum of
all tuples in the product. The component corresponding to $I\in
\N^m$ we denote by $P_K(I)$. The other is an $\N$-grading which
corresponds to the number of multipliers in the tensor product. We
will denote the components as $P_K=\oplus_{n\in \N}(P_K)_n$. The
total grading is calculated as $|I|-n$.

Notice that this construction is similar to the construction of
$T_K$, and so the inclusion $\{0,1\}\subset \N$ induces the
natural monomorphism $T_K\hookrightarrow P_K,$ so we will regard
$T_K$ as the subalgebra of $P_K$.

\begin{lemma}
The $\dg$-algebra $\widetilde{W}_K(\AAA)$ is isomorphic to the
following $\dg$-algebra:
$$G_K=\oplus_{I\in \N^m}P_K(I)\otimes (B_*{\AAA})^{I},$$
with the differential $d$ is the sum of 3 differentials
$d=d_{P}+d_{bar}+d_A$, where
\begin{itemize}
\item
$d_P$ corresponds to the differential in $P_K$: $$d_P(p\otimes
{\bf b})=d(p)\otimes {\bf b}.$$
\item
$d_{bar}$ is the differential coming from external bar
differentials on $B_*(A_i)$;
\item
$d_A$ corresponds to the internal differentials in the
$\dg$-algebras $A_i$'s:
$$d_A(p\otimes {\bf b})=\sum_{j,k}(-1)^{\epsilon_{j,k}}p\otimes d_{j,k}({\bf b}),$$
where ${\bf b}$ is regarded as a set of $b_{jk}\in A_j$ for $j\in
[m]$, $k\in [i_j]$, and $d_{j,k}$ is acting only on $b_{jk}\in
A_j$ as the internal differential in $A_j$.

\end{itemize}
\end{lemma}

Defining this isomorphism we should be careful with signs as they
depend on the gradings in ${\bf A}$ as well. We skip in the text the careful examining of those signs
just giving the resulting formulas.

$G_K$ is graded by usual rule for tensor products. We again split
this grading as follows. $G_K$ inherits multi-grading in $\N^m$
from $P_K$. The second grading of the element $p\otimes {\bf a}$
for $p\in (P_K)_n$ we define as $s=-n+\deg{\bf a}$ and refer to it
as $s$-grading. Respect to the last grading the differentials
satisfies the following conditions:
\begin{itemize}
\item
$d_{bar}$ preserves $s$-grading;
\item
$d_P$ and $d_A$ decrease the $s$-grading by 1.
\end{itemize}

Sometimes we will write the elements of the canonical basis of
$P_K$ in the form $(p(1)|\dots|p(s))$, where each $p(k)$ is a
monomial on $m$ variables $x_1,\dots,x_m$. We will consider any
element from $P_K$ as a linear combination of the basis elements,
and refer to them as summands, and to components $p(k)$ of these
summands as monomials.

In this terms, the subalgebra $T_K$ defined by the following
condition: it is spanned by all such elements of the canonical
basis such that all their monomials are square-free. We denote
this subalgebra also by $\PKN=T_K$.

Consider another condition on the powers of monomials: we take
$(p(1)|\dots|p(s))$ if all its monomials are square-free except
probably one of them (say, $p(k)$), and this exceptional monomial
$p(k)$ has degree 2 respect to one of the variables, and is
square-free respect to all the others (so it is of the form $p(k)=
x_{j_1}\dots x_{j_i}^2\dots x_{j_s}$). The vector subspace spanned
by such basis elements and $\PKN$ we denote by $P_K^1$. Note that
it is not a subalgebra of $P_K$.

Define the following operation $\nu_{j,k}$ on $P_K$: up to sign for each pair
$(j,k)$ where $j\in [m]$ and $k\in[i_j]$ we replace the $k$-th
entrance of $x_j$ by $x_j^2$. More precisely, let $\alpha=(p_1|x_j\tau|p_2)\in (P_K)_n$, where $p_1\in P_K(I)$ with $i_j=k-1$, and $p_2\in P_K$.
Then
$$\nu_{j,k}(\alpha)=(-1)^{\deg p_1+\epsilon(j,\tau)+n}(p_1|x_j^2\tau|p_2).$$

The definition of $P_K^1$ is
equivalent to the following: it is spanned by $\PKN$ and the
summands of the form $\nu_{j,k}(\alpha)$ where $\alpha\in \PKN$.

Thus, we constructed the inclusions
$$T_K=\PKN\subset P_K^1\subset P_K.$$
They induce the inclusions $$G_K^{\, 0}\subset G^1_K\subset G_K,$$
and $G_K^{\, 0}$ is a $\dg$-subalgebra of $G_K$.

Now the statement of Theorem~\ref{T:main} follows from the
following lemma which will be proved in the next subsections.

\begin{lemma}\label{Lem:epiG}
The homomorphisms $H_*(G^{\, 0}_K)\to H_*(G_K)$ and
$H_*(G^1_K/G^{\, 0}_K)\to H_*(G_K/G_K^{\, 0})$ induced by the
inclusions are surjective.
\end{lemma}

\begin{proof}[Proof of the Theorem~\ref{T:main}]
Denote the sequence of $\dg$-algebras ${\bf A}=C_*(\Omega \XX)$.

The first epimorphism of Lemma~\ref{Lem:epiG} any homology class
in $G_K$ has a representative in $G_K^{\, 0}$. For any element
$g_0\in \PKN(I)\otimes{\bf A}^I$ we have $d_{bar}=0$, and so the differential splits
respect to the splitting $G_K^{\, 0}\cong \PKN\otimes \AAA$.
Hence, the K\"unneth isomorphism implies
the required splitting of homology classes.

Now we should check that the equivalence relation can be defined
by $\mu_{j,k}$-operations.

We have that the differential restricted to $G^1_K/G^{\, 0}_K$ again splits,
and so due to the K\"unneth isomorphism it is enough to consider images under the differential of
elements of the form
$g=\nu_{j,k}(p)\otimes {\bf a}$ for some $p\in [\PKN(I)]_n$
($j\in[m]$, $k\in[i_j]$) and ${\bf a}=(\dots,a_j^k,
a_j^{k+1},\dots)\in {\bf A}^{I+e_j}$ with $d_P(p)=0$ and $d{\bf a}=0$.

Now we have $$dg=(-1)^n\left( \mu_{i,j}(p)\otimes {\bf a}-\nu_{j,k}p\otimes
{(\dots,a_j^k\cdot a_j^{k+1},\dots)}\right).$$
This gives the required
equivalence relations and finishes the proof of
Theorem~\ref{T:main}.
\end{proof}

The proof of Lemma~\ref{Lem:epiG} is proved in the rest of this
section. The first step is to reduce the statement about
surjectivity to the analogous question for the corresponding
homomorphisms in $P_K$ (Lemma~\ref{Lem:epi2}). The second step is
to show that the surjectivity for the homomorphism
$H_*(P^1_K/\PKN)\to H_*(P_K/\PKN)$ follows from the surjectivity
of the homomorphism $H_*(\PKN)\to H_*(P_K)$
(Lemma~\ref{Lem:epi1}). And finally we prove that the last
homomorphism is an epimorphism geometrically passing to diagonal
arrangements (Lemma~\ref{Lem:epi2}).

\subsection{Reduction to $P_K$}
The aim of this subsection is to prove the implications: the
surjectivity of the homomorphism $H_*(T_K)\twoheadrightarrow
H_*(P_K)$ (induced by the inclusion $T_K=\PKN\hookrightarrow P_K$)
implies the surjectivity of the homomorphism $H_*(G^{\, 0}_K)\to
H_*(G_K)$; and the surjectivity of the homomorphism
$H_*(P^1_K/\PKN)\twoheadrightarrow H_*(P_K/\PKN)$ implies the
surjectivity of the homomorphism $H_*(G^1_K/G^{\, 0}_K)\to
H_*(G_K/G^{\, 0}_K)$.

Both facts follow from the following technical statement.

Let $(G,d_G)$ be $\dg$-complex with an additional grading
$G=\oplus_{k\in \N} G_k$ and the differential split into sum of
two: $d_G=d_1+d_2$ in such a way that they satisfy the following
conditions respect to the additional grading:
\begin{enumerate}
\item
$d_1(G_k)\subset G_{k+1}$;
\item
$d_2(G_k)\subset G_k$.
\end{enumerate}

Further, suppose that $(G, d_G)$ is constructed as $G=\oplus_{I\in
\N^m} P(I)\otimes B(I)$ for some differential complexes
$(P=\oplus_{I\in\N^m} P(I), d_P)$ and $(B=\oplus_{I\in\N^m} B(I),
d_B)$, and the differential $d_1$ also splits:
$$d_1(p\otimes b)=d_P(p)\otimes b+(-1)^{\deg p}p\otimes d_B(b).$$

\begin{lemma}
Suppose we are given an embedding $\dg$-subalgebra
$P^0\hookrightarrow P$ such that (a) it induces epimorphism in
homology; (b) the restriction of $d_2$ on $P^0$ is zero. Then the
induced embedding of $\dg$-algebras $G^{\, 0}\hookrightarrow G$
also induces epimorphism in homology.
\end{lemma}
\begin{proof}
Given an element $g\in G$ with $dg=0$, decompose it into the sum
$g=\sum_{k\in\N} g_k$ respect to the additional grading, and find
the maximal integer $k$ such that the $k$-th homogenious component
$g_k$ of $g$ is not in $G^{\, 0}$.

As $g_{k+1}\in G^{\, 0}$ we have $d_2g_{k+1}=0$. Together with
$dg=0$ this implies $d_1g_k=0$. As the differential $d_1$ splits,
we can find such an element $g'\in G_{k-1}$ that $g_k+dg'=\sum
p'_i\otimes b'_i$ with $d_P(p'_i)=0$ and $d_B(b'_i)=0$ for each
$i$.

Now, by assumption of the lemma, for each $i$ there exists
$p''_i\in P_{k-1}$ such that
$$p_i'+d_P(p''_i)\in P^0.$$
Consider the element $\tilde{g}=g+d_G(g'+\sum_{i}p''_i\otimes
b'_i)\in G$. The element $\tilde{g}$ has the same components in
degrees $>k$ as $g$ has, but now we have $\tilde{g}_k\in G^{\,
0}$. So we decreased the value of $k$ and proceeding in this way
we get the required statement.
\end{proof}

\subsection{Lemmas~\ref{Lem:epi1} and \ref{Lem:epi2}}
In this subsection we prove that the embeddings
$T_K=\PKN\hookrightarrow P_K$ and $P_K^1/\PKN\hookrightarrow
P_K/\PKN$ induce the epimorphisms in homology.

\begin{lemma}\label{Lem:epi1}
The surjectivity of the homomorphism $H_*(\PKN)\to H_*(P_K)$
implies the surjectivity of the homomorphism $H_*(P^1_K/\PKN)\to
H_*(P_K/\PKN)$.
\end{lemma}

\begin{proof}
We need to prove that if for $\alpha\in P_K$ the condition
$d\alpha\in \PKN$ holds, then there exists $\beta\in P_K$ such
that $\alpha+d\beta\in P_K^1$.

Let  $\alpha\in P_K(I)$. We will prove the statement using the
induction by $(\sum_{s=1}^m i_s - m)$. The base of the induction
is trivial as when $I\preceq (1,\dots,1)$ the complexes $P_K(I)$
and $\PKN(I)$ coincide and so in this case $\alpha\in \PKN$.

Fix the number $j\in [m]$ and let $n=i_j$. So in each summand in
$\alpha$ we have exactly $n$ entries of $x_j$, let enumerate them
from left to right by $x_j^1$,\dots,$x_j^n$.

First, suppose that there exist an integer $k$ ($k\in [n]$) such
that there are no monomials in $\alpha$ divisible  by the product
$x_j^k x_j^{k+1}$ (let call this assumption a separating condition
respect to $x_j^k$ and $x_j^{k+1}$). Then we can rename fist $k$
entries of $x_j$ to $x_j'$, and next $n-k$ entries of $x_j$ to
$x_j''$. The resulting chain $\alpha'$ will be in $T_{K'}$ where
$K'=K(\pt,\dots,\pt,V[2],\pt,\dots,\pt)$ (here $V[2]$ is a
disjoint union of vertices $j'$ and $j''$, and it is the $j$-th
argument in the functor; the functor $K$ is defined in
\ref{Cons:functorK}). The differential of $\alpha'$ satisfies the
condition $d\alpha'\in P^{\,0}_{K'}$. Now by induction assumption
(for $\alpha'$ the number $\sum_{s=1}^m i_s - m$ is less by 1 than
for $\alpha$)
$$\alpha'+d\beta'\in P^{\, 0}_{K'}.$$
Renaming back $j'$ and $j''$ to $j$ we get the required statement
for $\alpha$.

Now consider the case when the separation assumption doesn't hold.
So fix some $k$ and renumber the index set
$\{1,\dots,k,k+1,k+2,\dots,n\}$ to
$\{1,\dots,k,k',k+1,\dots,n-1\}$.

Let split the complex $P_K$ as $F_1\oplus F_2$, where $F_1$ is
additively generated by the basis elements  where $x_j^k$ and
$x_j^{k'}$ are not separated, and $F_2$ --- by the basis elements
satisfying separating condition respect to $j$ and $k$. There is a
canonical isomorphism as additive groups
$P_K(\hat{I}):=P_K(i_1,\dots,n-1,\dots,i_m)\rightarrow F_1$,
denoted as before by $\nu_{j,k}$, which replace $x_j^k$ by the
product $x_j^k x_j^{k'}$ with certain sign. The following holds for the
differential:
$$d\,\nu_{j,k}(\lambda)=\nu_{j,k}(d\lambda) +(-1)^\delta \mu_{j,k}(\lambda),$$
where $\nu_{j,k}(d\lambda)\in F_1$,  $\mu_{j,k}(\lambda)\in F_2
\subset P_K(I)$ is defined similarly to $\mu_{j,k}^T$ from
Section~\ref{Subsec:doubling}, and the sign $(-1)^\delta$ will be not important here.

Let $\alpha=\alpha_1+\alpha_2$ be the decomposition corresponding
to the splitting constructed above. Then we have
$\alpha_1=\nu_{j,k}(\lambda)$ for some $\lambda\in
P_K(i_1,\dots,n-1,\dots,i_m)$. The condition $d\alpha=0$ implies
$d\lambda=0$. By the forthcoming Lemma~\ref{Lem:epi2} we have
$\lambda + d\Lambda=\gamma\in \PKN$, and, hence,
$$\nu_{j,k}({\lambda})+d\nu_{j,k}({\Lambda})=\nu_{j,k}({\gamma})+(-1)^\delta \mu_{j,k}(\Lambda).$$
We know that $d\nu_{j,k}({\gamma})\in \PKN$ and we have
$$\alpha + d\nu_{j,k}({\Lambda}) = \nu_{j,k}(\gamma) + (-1)^\delta \mu_{j,k}(\Lambda) + \alpha_2,$$
where (a) $\nu_{j,k}({\gamma})\in P_K^1$; (b)
$d\nu_{j,k}(\gamma)\in \PKN$; (c) $\alpha_2+(-1)^\delta \mu_{j,k}(\Lambda)\in
F_2$.

From (b) and the condition that $d\alpha\in \PKN$ it follows that
$d(\alpha_2+\nu_{j,k}(\Lambda))\in \PKN$. So as now
$\alpha_2+(-1)^\delta \mu_{j,k}(\Lambda)$ satisfies the condition of the lemma
and the separating assumption, by the first part of the proof
there exists $\beta$ such that $\alpha_2+(-1)^\delta \mu_{j,k}(\Lambda) +
d\beta\in P_K^1$ and so
$$\alpha+d\nu_{j,k}(\Lambda)+d\beta\in P_K^1.$$
This completes the proof.
\end{proof}

\begin{lemma}\label{Lem:epi2}
The homomorphism $H_*(\PKN)\to H_*(P_K)$ is surjective.
\end{lemma}
\begin{proof}

Fix the multi-degree $I=\{i_1,\dots,i_m\}$. Denote
$\tilde{K}=K(\Delta[i_1],\dots,\Delta[i_m])$ and
$K'=K(V[i_1],\dots,V[i_m])$. we will use the following notation: $
T_0=T_{K'}(1^{i_1},\dots,1^{i_m})$ and
$T=T_{\tilde{K}}(1^{i_1},\dots,1^{i_m})$. The inclusion of the
simplicial complexes  $K'\subset \tilde{K}$ induces the inclusion
of the differential graded chain complexes: $T_0\hookrightarrow
T$.

Consider the pair $(P_K(I), \PKN(I))$ and show that it is the
retract of the pair $(T,T_0)$. For that we construct the maps of
differential chain complexes $Pol: (P_K(I),
P^0_K(I)\hookrightarrow (T, T_0)$, and $Q: (T,T_0)\to
(P_K(I),\PKN(I))$ such that $Q\circ Pol=\id_{(P_K(I), \PKN(I))}$.

Let $\alpha$ be a summand in $P_K(I)$. First define the inclusion
$\varphi: P_K(I)\to T$ by the following procedure. For each $j$ we
have exactly $i_j$ entries of $x_j$ in $\alpha$, rename them from
left to right to $x_j^1$,\dots,$x_j^n$, considering new variables
as distinct. This map doesn't commute with the differential.

Now consider all {\it different} summands which are obtained from
$f(\alpha)$ by permuting the upper indices of $x_j$'s, and define
$Pol(\alpha)$ as the sum of them with the appropriate signs. Then
we have that $d Pol(\alpha) = Pol(d\alpha)$, and
$Pol(P^0_K(I))\subset T_0$.

Introduce the operation $Q'$ : it is identical on $f(P_K(I))$, and
it maps all the summands not from $P_K(I)$ to zero. Set
$Q=f^{-1}\circ Q$. We have $d Q(\beta)=Q(d\beta)$, and
$Q(T_0)=\PKN$.

Obviously, $Q\circ P=f^{-1}\circ Q'\circ P=f^{-1}\circ Q'\circ f=
f^{-1}\circ f=\id$. So we constructed the retraction.

Now recall that for any $L$ the differential complex
$T_L(1,\dots,1)$ is quasi-isomorphic to the chain complex of $D_L$
- the complement of  the corresponding diagonal arrangement, and
for $L_1\subset L_2$ on the same set of vertices the inclusion
$T_{L_1}(1,\dots,1)\hookrightarrow T_{L_2}(1,\dots,1)$ corresponds
to the inclusion $D_{L_1}\hookrightarrow D_{L_2}$. So by the next
Lemma~\ref{Lem:diag1} we have that $T_0\hookrightarrow T$ induces
an epimorphism in homology. As $(P_K(I),\PKN(I))$ is the retract
of the pair $(T,T_0)$, the same is true for the homomorphism
$H_*(P^0_K(I))\to H_*(P_K(I))$ - it is surjective.

\end{proof}

\begin{lemma}\label{Lem:diag1}
For $\tilde{K}=K(\Delta[i_1],\dots,\Delta[i_m])$ and
$K'=K(V[i_1],\dots,V[i_m])$ the inclusion $D_{K'}\hookrightarrow
D_{\tilde{K}}$ induces an epimorphism in homology.
\end{lemma}
\begin{proof}
Inductive arguments shows that it is enough to prove the following
lemma.

\begin{lemma}\label{Lem:edge}
Let $\{i,j\}$ be an edge in $L$ and let $\hat{L}$ denotes the
complex obtained from $L$ by deleting all simplices containing the
 edge $\{i,j\}$. Let $L$ satisfy the following condition: for
any $\sigma\in K$ such that $i\in \sigma$, the subset
$\sigma\cup\{j\}$ is also in $K$. Then $D_{\hat{L}}\hookrightarrow
D_L$ induces an epimorphism in homology.
\end{lemma}

Indeed, using this statement we can pass from any simplicial
complex of the form
$K_1=K(L_1,\dots,V[S\cup\{i\}]\star\Delta[T\cup\{j\}],L_m)$ to the
simplicial complex
$K_1=K(L_1,\dots,V[S\cup\{i,j\}]\star\Delta[T],L_m)$, $\star$
denotes the join of simplicial complexes.

So let now prove Lemma~\ref{Lem:edge}.

Let $\alpha$ is a $k$-dimensional cycle in $D_L$. We can assume
that it is transversal to the hyperplane $\pi$ defined by the
equation $x_i=x_j$. Let $\alpha_+$ denotes the intersection of
$\alpha$ with the halfspace $x_i\ge x_j$, and $\alpha_-$ with the
halfspace $x_i\le x_j$. The intersection $\alpha\cap \pi$ defines
a $(k-1)$-dimensional cycle $\lambda$ in $D_L\cap \pi\cong
D_{\hat{L}}$ with $d\alpha_+=-d\alpha_-=i(\lambda)$, where $i:
D_L\cap \pi \to D_L$ is a natural inclusion.

Consider the projection $P: D_L\rightarrow D_{\hat{L}}$ which is
forgetting the coordinate $x_j$. Denote $\beta=P_*(\alpha_+)$. As
$P_*(i(\lambda))=\lambda$ we have $d\beta=\lambda$. So $d(\alpha_+
+ i(\beta)) = d(\alpha_--i(\beta))=0$. We have the following
decomposition of $\alpha$ into the sum of two cycles
$$\alpha=(\alpha_++i(\beta)) + (\alpha_- - i(\beta)).$$
But adding a small $\varepsilon$ to the $x_i$-coordinate of all
points of the first cycle in the sum, and subtracting also
$\varepsilon$ to the $x_i$-coordinates of all points of the second
cycle we deform them to homology equivalent cycles which are
contained in open halfspaces defined by inequalities $x_i>x_j$ and
$x_i<x_j$ respectively.
\end{proof}


\begin{thebibliography}{100}
\bibitem{anick}
Anick D. Connections between Yoneda and Pontrjagin algebras.
Algebraic topology, Aarhus 1982 (Aarhus, 1982), 331--350, Lecture
Notes in Math., 1051, Springer, Berlin, 1984.
\bibitem{cohenetc}
Bahri A., Bendersky M., Cohen F. R., and Gitler S.
The polyhedral product functor: a method of computation for moment-angle complexes, arrangements and related spaces, 2007, arXiv: 0711.4689.

\bibitem{baryshnikov}
Baryshnikov Y. Cohomology ring of no k-equal manifolds, preprint,
1997.
\bibitem{bjorner2}
Bj\"orner A. Subspace arrangements. First European Congress of
Mathematics, Vol. I (Paris, 1992), 321--370, Progr. Math., 119,
Birkhauser, Basel, 1994.
\bibitem{bjorner}
Bj\"orner A., Welker V.  The homology of ``$k$-equal'' manifolds
and related partition lattices. Adv. Math. 110 (1995), no. 2,
277--313.
\bibitem{bod}
B\"odigheimer C.-F. Stable splittings of mapping spaces. Lecture
Notes in Math., 1286, Springer, Berlin, 1987 , 174--187,
\bibitem{BP:book}
Buchstaber V., Panov T. Torus Actions and Their Applications in
Topology and Combinatorics, AMS University Lectures Series,
vol.24, 2002.
\bibitem{CS}
Clark A., Smith L. The rational homotopy of a wedge. Pacific J.
Math. 24 1968 241--246.
\bibitem{GT}
Grbic J., Theriault S. The homotopy type of the complement of a
coordinate subspace arrangement. Topology 46 (2007), no. 4,
357--396.
\bibitem{dobrf}
Dobrinskaya, N.
Configuration spaces of labelled particles and finite Eilenberg-MacLane complexes.
Proc. Steklov Inst. Math. 2006, no. 1 (252), 30--46
\bibitem{Dobr:Poisson}
Dobrinskaya N. Generalized Poisson structures in loop space
homology. In preparation.
\bibitem{lemaire}
Lemaire, J.-M.
Alg\`ebres connexes et homologie des espaces de lacets.
Lecture Notes in Mathematics, Vol. 422, 1974.
\bibitem{NR}
Notbohm D., Ray N. On Davis-Januszkiewicz homotopy types. I.
Formality and rationalisation. Algebr. Geom. Topol. 5 (2005),
31--51
\bibitem{panray}
Panov T., Ray R. Categorical aspects of toric topology.
Contemporary Mathematics series, 460 (2008), 298 -322
\bibitem{papasuc}
Papadima S, Suciu A.
Algebraic invariants for right-angled Artin groups.
Math. Ann. 334 (2006), no. 3, 533--555.
\bibitem{porter}
Porter G. The homotopy groups of wedges of suspensions. Amer. J.
Math. 88, 1966, 655--663.
\bibitem{segal}
Segal G. Configuration-spaces and iterated loop-spaces. Invent.
Math. 21 (1973), 213--221.
\bibitem{PRW}
Peeva I., Reiner V., Welker V.
Cohomology of real diagonal subspace arrangements via resolutions.
Compositio Math. 117 (1999), no. 1, 99--115.
\end{thebibliography}
\end{document}